\documentclass[12pt,reqno]{amsart}
\usepackage{amsmath,amsthm,amsfonts}
\usepackage[mathscr]{eucal}
\usepackage{eufrak}
\usepackage{graphicx}

{\obeylines
\gdef\MATH{\begingroup\parindent0pt\parskip0pt plus 0pt\obeylines%
        \def^^M{\vskip4pt}%
        \obeyspaces\tt\small}%
}
\def\goodbreakpoint{\par\penalty-5000%
         \vrule height10pt depth2pt width0pt\leavevmode}
\def\endMATH{\endgroup}
\def\MATHphi{\leavevmode
        \hbox to 0pt{\hbox to 5.24995pt{\hss$\phi$\hss}\hss}}
\def\MATHGamma{\leavevmode
        \hbox to 0pt{\hbox to 5.24995pt{\hss$\Gamma$\hss}\hss}}
\def\MATHpi{\leavevmode
        \hbox to 0pt{\hbox to 5.24995pt{\hss$\pi$\hss}\hss}}
\def\MATHinfty{\leavevmode
        \hbox to 0pt{\hbox to 5.24995pt{\hss$\infty$\hss}\hss}}
\def\MATHhStrich{\leavevmode
        \hbox to 0pt{\hbox to 5.24995pt{\vrule height4.5pt depth-3.5pt width5.24995pt}\hss}}
\def\MATHluEck{\leavevmode
        \hbox to 0pt{\hbox to 5.24995pt{\hskip2.12497pt
         \vrule height4.5pt depth1pt width1pt
         \vrule height4.5pt depth-3.5pt width2.12498pt}\hss}}
\def\MATHruEck{\leavevmode
        \hbox to 0pt{\hbox to 5.24995pt{%
         \vrule height4.5pt depth-3.5pt width2.12497pt
         \vrule height4.5pt depth1pt width1pt
         \hskip2.12498pt}\hss}}
\def\MATHloEck{\leavevmode
        \hbox to 0pt{\hbox to 5.24995pt{\hskip2.12497pt
         \vrule height9pt depth-3.5pt width1pt
         \vrule height4.5pt depth-3.5pt width2.12498pt}\hss}}
\def\MATHroEck{\leavevmode
        \hbox to 0pt{\hbox to 5.24995pt{%
         \vrule height4.5pt depth-3.5pt width2.12497pt
         \vrule height9pt depth-3.5pt width1pt
         \hskip2.12498pt}\hss}}
\def\MATHvStrich{\leavevmode
        \hbox to 0pt{\hbox to 5.24995pt{\hskip2.12497pt
         \vtop to 0pt{\hsize1pt\vss%
                \vrule height17pt depth6pt width1pt\vskip8pt\vss\par}%
         \hskip2.12498pt}\hss}}
\def\MATHtStueck{\leavevmode
        \hbox to 0pt{\hbox to 5.24995pt{%
         \vrule height4.5pt depth-3.5pt width2.12497pt
         \vrule height4.5pt depth2pt width1pt
         \vrule height4.5pt depth-3.5pt width2.12498pt}\hss}}
\def\MATHbackslash{\leavevmode
        \hbox to 0pt{\hbox to 5.24995pt{\hss$\backslash$\hss}\hss}}
\def\MATHlbrace{\leavevmode
        \hbox to 0pt{\hbox to 5.24995pt{\hss$\{$\hss}\hss}}
\def\MATHrbrace{\leavevmode
        \hbox to 0pt{\hbox to 5.24995pt{\hss$\}$\hss}\hss}}
\def\MATHkleiner{\leavevmode
        \hbox to 0pt{\hbox to 5.24995pt{\hss$\langle$\hss}\hss}}
\def\MATHkleiner{<}
\def\MATHgroesser{\leavevmode
        \hbox to 0pt{\hbox to 5.24995pt{\hss$\rangle$\hss}\hss}}
\def\MATHgroesser{>}
\def\MATHhoch{\leavevmode
        \hbox to 0pt{\hbox to 5.24995pt{\hss$^\land$\hss}\hss}}
\def\MATHtief{\leavevmode
        \hbox to 0pt{\hbox to 5.24995pt{\hss\vrule height0pt depth.8pt width3pt\hss}\hss}}

\allowdisplaybreaks
\newtheorem{theorem}{Theorem}
\theoremstyle{plain}

\newtheorem{lemma}{Lemma}

\numberwithin{equation}{section}
\input{tcilatex}

\begin{document}
\title[An eigenvalue problem]{An eigenvalue problem for\\
the associated Askey--Wilson polynomials}
\author{Andrea Bruder}
\address{Mathematics and Computer Science, Colorado College, Tutt Science
Center, 14 E. Cache la Poudre St., Colorado Springs, CO 80903, U.S.A.}
\email{Andrea.Bruder@coloradocollege.edu}
\author{Christian Krattenthaler$^\dagger$}
\address{Fakult\"{a}t f\"{u}r Mathematik, Universit\"{a}t Wien, Nordbergstra%
\ss e 15, A-1090 Vienna, Austria}
\urladdr{http://www.mat.univie.ac.at/\lower0.5ex\hbox{\~{}}kratt}
\author{Sergei K. Suslov}
\address{School of Mathematical and Statistics Sciences \& Mathematical.
Computational and Modeling Sciences Center, Arizona State University, Tempe,
AZ 85287-1804, U.S.A.}
\email{sks@asu.edu}
\urladdr{http://hahn.la.asu.edu/\symbol{126}suslov/index.html}
\date{\today }
\thanks{$^\dagger$Research partially supported by the Austrian Science
Foundation FWF, grants Z130-N13 and S9607-N13, the latter in the framework
of the National Research Network ``Analytic Combinatorics and Probabilistic
Number Theory"}
\subjclass{33D45, 42C10; 33D15}
\keywords{Basic hypergeometric functions, $q$-orthogonal polynomials,
Askey--Wilson polynomials, associated Askey--Wilson polynomials, eigenvalue
problem.}

\begin{abstract}
To derive an eigenvalue problem for the associated Askey--Wilson
polynomials, we consider an auxiliary function in two variables which is
related to the associated Askey--Wilson polynomials introduced by Ismail and
Rahman. The Askey--Wilson operator, applied in each variable separately,
maps this function to the ordinary Askey--Wilson polynomials with different
sets of parameters. A third Askey--Wilson operator is found with the help of
a computer algebra program which links the two, and an eigenvalue problem is
stated.
\end{abstract}

\maketitle

\section{Introduction}

Throughout this paper, we use the standard notation for the $q$-shifted
factorials:
\begin{alignat*}2
\left( a;q\right) _{n}&:=\prod\limits_{j=0}^{n-1}\left( 1-aq^{j}\right) ,&%
\qquad \left( a_{1},a_{2},\dots,a_{r};q\right)
_{n}&:=\prod_{k=1}^{r}\left( a_{k};q\right) _{n},  %\label{in1}
\\
\left( a;q\right) _{\infty }&:=\lim_{n\rightarrow \infty }\left( a;q\right)
_{n},&\qquad \left( a_{1},a_{2},\dots,a_{r};q\right) _{\infty
}&:=\prod_{k=1}^{r}\left( a_{k};q\right) _{\infty },  %\label{in2}
\end{alignat*}%
provided $\left\vert q\right\vert <1.$ The basic hypergeometric series is
defined by (cf.\ \cite{Ga:Ra})
\begin{equation*}
_{r}\varphi _{s}\left(
\begin{array}{c}
a_{1},a_{2},\dots,a_{r} \\
b_{1},\dots,b_{s}%
\end{array}%
;\,q\,,\,z\right) :=\sum_{n=0}^{\infty }\frac{\left(
a_{1},a_{2},\dots,a_{r};q\right) _{n}}{\left(
q,b_{1},b_{2},\dots,b_{s};q\right) _{n}}\,((-1)^{n}q^{n(n-1)/2})^{1+s-r}%
\,z^{n}. 
\end{equation*}%
If $0<|q|<1,$ the series converges absolutely for all $z$ if $r\leq s,$ and
for $|z|<1$ if $r=s+1.$

The Askey--Wilson polynomials are the most general extension of the classical
orthogonal polynomials \cite{An:As}, \cite{An:As:Ro}, \cite{As:Wi}, \cite%
{Koe:Sw}, \cite{Ni:Su:Uv}, \cite{Sz}. They are most conveniently
given in terms of a $_4\varphi_3$-series,
\begin{align} 
p_{n}(x)& =p_{n}(x;a,b,c,d)=p_{n}(x;a,b,c,d|q)   \notag \\
& =a^{-n}\,(ab,ac,ad;q)_{n}\;{}_{4}\varphi _{3}\!\left( \!\!%
\begin{array}{c}
q^{-n},\ abcdq^{n-1},\ az,\ a/z\smallskip \\[0.1cm]
ab,\ ac,\ ad%
\end{array}%
\!\!;q,\,q\!\right) ,  \notag
\end{align}%
where $x=\left( z+z^{-1}\right) /2,$ and $\left\vert z\right\vert <1.$ In this
normalization, the Askey--Wilson polynomials are symmetric in all four
parameters due to Sears' transformation \cite{As:Wi}.

The Askey--Wilson polynomials satisfy the 3-term recurrence relation%
\begin{equation}
2x\,\,p_{n}(x;a,b,c,d)=A_{n}\,p_{n+1}(x;a,b,c,d)\,+\,B_{n}%
\,p_{n}(x;a,b,c,d)+\,C_{n}\,p_{n-1}(x;a,b,c,d),  \label{in5}
\end{equation}%
where%
\begin{align}
A_{n}& =\frac{a^{-1}(1-abq^{n})(1-acq^{n})(1-adq^{n})(1-abcdq^{n-1})}{%
(1-abcdq^{2n-1})(1-abcdq^{2n}-q^{2n})},  \label{in6} \\
C_{n}& =\frac{a(1-bcq^{n-1})(1-bdq^{n-1})(1-cdq^{n-1})(1-q^{n})}{%
(1-abcdq^{2n-1})(1-abcdq^{2n})},  \label{in7} \\
B_{n}& =a+a^{-1}-A_{n}-C_{n}.  \label{in8}
\end{align}%
The weight function with respect to which the polynomials $p_{n}(x)$ are
orthogonal was found by Askey and Wilson in \cite{As:Wi}. The Askey--Wilson
divided difference operator is defined by%
\begin{align}
L(x)u&:=L\left( s;a,b,c,d\right) u\left( s\right)   \notag \\
&\hphantom{:} =\frac{\sigma \left( -s\right) \nabla x\left( s\right) u\left(
s+1\right) +\sigma \left( s\right) \Delta x\left( s\right) u\left(
s-1\right) -\left[ \sigma \left( s\right) \Delta x\left( s\right) +\sigma
\left( -s\right) \nabla x\left( s\right) \right] u\left( s\right) }{\Delta
x\left( s\right) \nabla x\left( s\right) \nabla x_{1}\left( s\right) },
\label{in9}
\end{align}%
where $\sigma \left( s\right) =q^{-2s}\left( q^{s}-a\right) \left(
q^{s}-b\right) \left( q^{s}-c\right) \left( q^{s}-d\right) $ and, by
definition,
\begin{align*}
x(s)& =\frac{1}{2}\left( q^{s}+q^{-s}\right) \text{\ }, & \qquad x_{1}(s)&
=x\left( s+\frac{1}{2}\right) , \\
\Delta f(s)& =f(s+1)-f(s), & \qquad \nabla f(s)& =f(s)-f(s-1).
\end{align*}%
(We follow the notation in \cite{At:Su:DHF} and \cite{At:Su1}.) We will make
use of an analogue of the power series expansion method, where a function is
expanded in terms of generalized powers. For a positive integer $m,$ the
generalized powers are defined by%
\begin{equation}
\lbrack
x(s)-x(z)]^{(m)}=\prod_{n=0}^{m-1}[x_{n}(s)-x_{n}(z-k)],\qquad
x_{n}(z)=x\left( z+\frac{n}{2}\right)   \label{in10}
\end{equation}%
(see \cite[Exercises~2.9--2.11, 2.25]{Su4} and \cite{Su2} for more
details).

\section{The Associated Askey--Wilson Polynomials}

The associated Askey--Wilson polynomials, $p_{n}^{\alpha }(x)=p_{n}^{\alpha
}(x;a,b,c,d)=p_{n}^{\alpha }(x;a,b,c,d|q),$ were introduced by Ismail and
Rahman in \cite{Is:Rah}. They are solutions of the 3-term recurrence relation%
\begin{equation}
2x\,\,p_{n}^{\alpha }(x;a,b,c,d)=A_{n+\alpha }\,\,p_{n+1}^{\alpha
}(x;a,b,c,d)\,+\,B_{n+\alpha }\,\,p_{n}^{\alpha }(x;a,b,c,d)+\,C_{n+\alpha
}\,\,p_{n-1}^{\alpha }(x;a,b,c,d),  \label{aaw1}
\end{equation}%
where $0<\alpha <1,$ with initial values \ $\,p_{-1}^{\alpha }(x)=0,$ $%
p_{0}^{\alpha }(x)=1$, and $A_{n+\alpha },$ $B_{n+\alpha },$ $C_{n+\alpha }$
are given as in (\ref{in6})--(\ref{in8}) with $n$ replaced by $n+\alpha .$
The two linearly independent solutions to (\ref{in5}) found in \cite%
{Is:Rah} are%
\begin{multline}
R_{n+\alpha } =\frac{(abq^{n+\alpha },acq^{n+\alpha },adq^{n+\alpha
},bcdq^{n+\alpha }/z;q)_{\infty }}{(bcq^{n+\alpha },bdq^{n+\alpha
},cdq^{n+\alpha },azdq^{n+\alpha };q)_{\infty }}\left( \frac{a}{z}\right)
^{n+\alpha }\smallskip \smallskip \medskip  \label{aaw2} \\
 \times \,_{8}W_{7}(bcd/qz;b/z,c/z,d/z,abcdq^{n+\alpha
-1},q^{-\alpha -n};q,qz/a)
\end{multline}%
and%
\begin{multline}
S_{n+\alpha } =\frac{(abcdq^{2n+2\alpha },bzq^{n+\alpha +1},czq^{n+\alpha
+1},dzq^{n+\alpha +1},bcdzq^{n+\alpha +1};q)_{\infty }}{(bcq^{n+\alpha
},bdq^{n+\alpha },cdq^{n+\alpha },q^{n+\alpha +1},bcdzq^{2n+2\alpha
+1};q)_{\infty }}(az)^{n+\alpha }\smallskip \smallskip \medskip  \label{aaw3}
\\
\times \,_{8}W_{7}(bcdzq^{2n+2\alpha };bcq^{n+\alpha
},bdq^{n+\alpha },cdq^{n+\alpha },q^{n+\alpha +1},zq/a;q,az).
\end{multline}%
The weight function for the associated Askey--Wilson polynomials and an
explicit polynomial representation were found by Ismail and Rahman in \cite%
{Is:Rah}. The latter is given by%
\begin{align}
p_{n}^{\alpha }(x)& =p_{n}^{\alpha }(x;a,b,c,d|q)  \notag \\
& =\sum_{k=0}^{n}\frac{(q^{-n},abcdq^{2\alpha +n-1},abcdq^{2\alpha
-1},ae^{i\theta },ae^{-i\theta };q)_{k}}{(q,abq^{\alpha },acq^{\alpha
},adq^{\alpha },abcdq^{\alpha -1};q)_{k}}\ q^{k}  \notag \\
& \qquad \times ~_{10}W_{9}(abcdq^{2\alpha +k-1};q^{\alpha },bcq^{\alpha
-1},bdq^{\alpha -1},cdq^{\alpha -1},q^{k+1},abcdq^{2\alpha
+n+k-1},q^{k-n};q,a^{2}).  \label{aaw4}
\end{align}%
There is another useful representation of the associated Askey--Wilson
polynomials in terms of a double series due to Rahman,
\begin{align}
p_{n}^{\alpha }(x)& =p_{n}^{\alpha }(x;a,b,c,d|q)\smallskip  \notag \\
& =\frac{(abcdq^{2\alpha -1},q^{\alpha +1};q)_{n}}{(q,abcdq^{\alpha
-1};q)_{n}}q^{-\alpha n}\sum_{k=0}^{n}\frac{(q^{-n},abcdq^{2\alpha
+n-1};q)_{k}}{(q^{\alpha +1},abq^{\alpha };q)_{k}}\smallskip  \notag \\
& \qquad \times \frac{(aq^{\alpha }e^{i\theta },aq^{\alpha }e^{-i\theta
};q)_{k}}{(acq^{\alpha },acq^{\alpha };q)_{k}}\sum_{j=0}^{k}\frac{(q^{\alpha
},abq^{\alpha -1},acq^{\alpha -1},adq^{\alpha -1};q)_{j}}{(q,abcdq^{2\alpha
-2},aq^{\alpha }e^{i\theta },aq^{\alpha }e^{-i\theta };q)_{j}}q^{j},
\label{aaw5}
\end{align}%
where $x=\cos \theta $ (see \cite[Exercises~8.26--8.27]{Ga:Ra} and \cite%
{Rahman96}, \cite{Rah2000}). This formula will be the starting point for our
investigation.

\section{An Overview of the Main Result}

To construct an eigenvalue problem for the associated Askey--Wilson
polynomials, let us consider an auxiliary function $u_{n}^{\alpha }(x,y)$ in
two variables, which for $x=y$ coincides with the associated Askey--Wilson
polynomials (up to a factor). We observe that the Askey--Wilson operator $%
L_{0}(x)$ (in one variable $x$) maps $u_{n}^{\alpha }(x,y)$ to the $n$-th
degree ordinary Askey--Wilson polynomial (up to some factors). A similar
result is obtained for the operator $L_{1}(y)$ applied to $u_{n}^{\alpha
}(x,y)$ with respect to the second independent variable $y.$ We will find an 
operator $L_{2}(x),$ which maps certain multiples of $\left(
L_{1}(y)+\lambda \right) u_{n}^{\alpha }(x,y)$\ to $(L_{0}(x)+\lambda
)u_{n}^{\alpha }(x,y).$ As a result, we obtain an eigenvalue problem of the
form%
\begin{multline}
\frac{(aq^{s},aq^{-s};q)_{\infty }}{(aq^{\alpha +s-1},aq^{\alpha
-s-1};q)_{\infty }}(L_{2}(x)+\lambda )\frac{(aq^{\alpha +s},aq^{\alpha
-s};q)_{\infty }}{(aq^{s},aq^{-s};q)_{\infty }}\left( L_{1}(y)+\mu _{\alpha
}\right) u_{n}^{\alpha }(x,y)  \label{ep1} \\
=\frac{4q^{9/2}}{(1-q)^{2}\gamma }(L_{0}(x)+\lambda _{\alpha
+n})u_{n}^{\alpha }(x,y)
\end{multline}%
related to the associated Askey--Wilson polynomials of Ismail and Rahman
(see Theorem~\ref{thm:1} below for an exact statement). We shall use the
normalization
\begin{equation}
p_{n}(x;a,b,c,d)={}_{4}\varphi _{3}\!\left( \!\!%
\begin{array}{c}
q^{-n},\ abcdq^{n-1},\ aq^{s},\ aq^{-s}\smallskip \\[0.1cm]
ab,\ ac,\ ad%
\end{array}%
\!\!;q,\,q\!\right)  \label{AskeyWilsonPlinomials}
\end{equation}%
for the ordinary Askey--Wilson polynomials
throughout this paper.

\begin{lemma} \label{lem:1}
Let $u_{n}^{\alpha }(x,y)$ be the function in the two variables $x$
and $y$ defined by
\begin{align}
u_{n}^{\alpha }(x,y)& :=\frac{(aq^{s},aq^{-s},aq^{\alpha +z},aq^{\alpha
-z};q)_{\infty }}{(aq^{\alpha +s},aq^{\alpha -s},aq^{z},aq^{-z};q)_{\infty }}
\notag \\
& \qquad \times \sum_{m=0}^{n}\frac{(q^{-n},\gamma q^{2\alpha
+n-1},aq^{\alpha +s},aq^{\alpha -s};q)_{m}}{(q^{\alpha +1},abq^{\alpha
},acq^{\alpha },adq^{\alpha };q)_{m}}q^{m}  \notag \\
& \qquad \qquad \times \sum_{k=0}^{m}\frac{(q^{\alpha },abq^{\alpha
-1},acq^{\alpha -1},adq^{\alpha -1};q)_{k}}{(q,\gamma q^{2\alpha
-2},aq^{\alpha +z},aq^{\alpha -z};q)_{k}}q^{k}, \label{AssAWFuncs}
\end{align}%
with $x(s)=(q^{s}+q^{-s})/2$ and $y(z)=(q^{z}+q^{-z})/2.$ Then $%
u_{n}^{\alpha }(x,y)$ satisfies an equation of the form
\begin{equation}
(L_{0}(x)+\lambda _{\alpha +n})u_{n}^{\alpha }(x,y)=f_{n}^{\alpha }(x,y),
\label{AWOpLemma1}
\end{equation}%
where $L_{0}(x)=L\left( s;a,b,c,d\right) $ is the Askey--Wilson divided
difference operator in the variable $x$ given by \eqref{in9}. Here,%
\begin{align}
f_{n}^{\alpha }(x,y)& =-\frac{4q^{3/2-\alpha }}{(1-q)^{2}}\frac{%
(aq^{s},aq^{-s},aq^{\alpha +z},aq^{\alpha -z};q)_{\infty }}{(aq^{\alpha
+s-1},aq^{\alpha -s-1},aq^{z},aq^{-z};q)_{\infty }}  \notag \\
& \qquad \quad \quad
\times (q^{\alpha },abq^{\alpha -1},acq^{\alpha -1},adq^{\alpha
-1};q)_{1}\,\,  \notag \\
& \qquad \qquad \times p_{n}(x;aq^{\alpha -1},bcdq^{\alpha
-1},q^{1+z},q^{1-z}),  \notag
\end{align}%
and%
\begin{equation*}
\lambda _{\alpha +n}=\frac{4q^{3/2}}{\left( 1-q\right) ^{2}}\left(
1-q^{-\alpha -n}\right) \left( 1-\gamma q^{\alpha +n-1}\right) ,\qquad
\gamma =abcd.
\end{equation*}
\end{lemma}

Note that $f_{n}^{\alpha }(x,y)$\ contains the $n$-th degree ordinary
Askey--Wilson polynomial of the form (\ref{AskeyWilsonPlinomials}) in the
variable $x.$ Our function $u_{n}^{\alpha }(x,y)$ is the Askey--Wilson
polynomial when $\alpha =0$ and a constant multiple of the associated
Askey--Wilson polynomial if $x=y.$

\begin{lemma} \label{lem:2}
The function $u_{n}^{\alpha }(x,y)$ satisfies another equation, namely%
\begin{equation*}
(L_{1}(y)+\mu _{\alpha })u_{n}^{\alpha }(x,y)=g_{n}^{\alpha }(x,y),
\end{equation*}%
where $L_{1}(y):=L\left( y;q/a,q/b,q/c,q/d\right) $ is the Askey--Wilson
divided difference operator in $y.$ Here,%
\begin{align}
g_{n}^{\alpha }(x,y)& =-\frac{4q^{9/2-\alpha }}{(1-q)^{2}\gamma }\frac{%
(aq^{s},aq^{-s},aq^{\alpha +z+1},aq^{\alpha -z+1};q)_{\infty }}{(aq^{\alpha
+s},aq^{\alpha -s},aq^{z},aq^{-z};q)_{\infty }}  \notag \\
& \qquad\quad \quad
\times (q^{\alpha },abq^{\alpha -1},acq^{\alpha -1},adq^{\alpha
-1};q)_{1}\,\,  \notag \\
& \qquad \qquad \times p_{n}(x;aq^{\alpha },bcdq^{\alpha -2},q^{1+z},q^{1-z})
\notag
\end{align}%
and%
\begin{equation*}
\mu _{\alpha }=\frac{4q^{3/2}}{\left( 1-q\right) ^{2}}\left( 1-q^{\alpha
}\right) \left( 1-q^{3-\alpha }/\gamma \right) .
\end{equation*}
\end{lemma}

Note that $g_{n}^{\alpha }(x,y)$\ contains another $n$-th degree
Askey--Wilson polynomial (\ref{AskeyWilsonPlinomials}) in the same variable $%
x.$

\begin{lemma} \label{lem:3}
The difference differentiation formula
\begin{equation}
(L\left( x\right) +\lambda )p_{n}(x;a,b,c,d)=\lambda p_{n}(x;a/q,bq,c,d)
\label{DiffDiffFormulaAW}
\end{equation}%
holds
for the Askey--Wilson polynomials given by \eqref{AskeyWilsonPlinomials}.
Here, $L\left( x\right) =L\left( s;a,a/q,c,d\right) $ is the Askey--Wilson
divided difference operator \eqref{in9} and%
\begin{equation*}
\lambda =\frac{4q^{3/2}}{\left( 1-q\right) ^{2}}\left( 1-ac/q\right) \left(
1-ad/q\right) .
\end{equation*}
\end{lemma}

Lemmas~\ref{lem:1}--\ref{lem:3}
allow us to establish the eigenvalue problem (\ref{ep1}) for the
associated Askey--Wilson functions (\ref{AssAWFuncs}),
see the next section.

\section{Main Result}
\label{sec:5}

With the help of Lemmas~\ref{lem:1}--\ref{lem:3},
we now identify an operator $L_{2}(x)$
linking $(L_{0}(x)+\lambda _{\alpha +n})u_{n}^{\alpha }(x,y)$ and $%
(L_{1}(y)+\lambda _{-\alpha })u_{n}^{\alpha }(x,y)$ in such a way that an
eigenvalue problem is formulated.

\begin{theorem} \label{thm:1}
Let $L_{2}(x)=L(s;aq^{\alpha },aq^{\alpha -1},q^{1+z},q^{1-z})$ be the
Askey--Wilson divided difference operator defined by \eqref{in9} with
\begin{equation*}
\sigma (s)=q^{-2s}\left( q^{s}-aq^{\alpha }\right) \left( q^{s}-aq^{\alpha
-1}\right) \left( q^{s}-q^{1+z}\right) \left( q^{s}-q^{1-z}\right)
\end{equation*}%
and%
\begin{equation*}
\lambda =\frac{4q^{3/2}}{(1-q)^{2}}\left( 1-aq^{\alpha -z}\right) \left(
1-aq^{\alpha +z}\right) .
\end{equation*}%
Then an eigenvalue problem for the associated Askey--Wilson functions $%
u_{n}^{\alpha }(x,y)$ can be stated as%
\begin{multline}
\frac{\gamma }{q^{3}}\frac{\left( aq^{s},aq^{-s};q\right) _{\infty }}{%
\left( aq^{\alpha +s-1},aq^{\alpha -s-1};q\right) _{\infty }}\left(
L_{2}(x)+\lambda \right) \frac{\left( aq^{\alpha +s},aq^{\alpha -s};q\right)
_{\infty }}{\left( aq^{s},aq^{-s};q\right) _{\infty }}\left( L_{1}(y)+\mu
_{\alpha }\right) u_{n}^{\alpha }(x,y)  \label{MainTh3} \\
=\frac{4q^{3/2}}{(1-q)^{2}}\left( L_{0}(x)+\lambda _{\alpha
+n}\right) u_{n}^{\alpha }(x,y),
\end{multline}%
where $L_{0},$ $L_{1},$ $\lambda _{\alpha +n},$ $\mu _{\alpha }$ and $%
u_{n}^{\alpha }(x,y)$ are defined as in Lemmas~\emph{\ref{lem:1}--\ref{lem:3}}.
\end{theorem}

Computational details are left to the reader. The explicit form of the
difference operator in two variables on the left-hand side of the last
equation has also been calculated, but it is too long to be displayed here.

\section{Proofs}

\subsection{Proof of Lemma \protect\ref{lem:1}}

Let $\lambda _{\nu }$ be an arbitrary number. We are looking for solutions
of a generalization of the equation~(\ref{AWOpLemma1}), namely,%
\begin{equation*}
(L_{0}(x)+\lambda _{\nu })u_{n}^{\alpha }(x,y)=f_{n}^{\alpha }(x,y),
\end{equation*}%
in terms of generalized powers (see \eqref{in10} for the definition)
\begin{equation*}
u_{n}^{\alpha }(x,y)=\sum_{m=0}^{n}c_{m}v_{m}[x(s)-x(\xi )]^{(\alpha +m)},
\end{equation*}%
where
\begin{equation*}
v_{m}=v_{m}(y)=\frac{(aq^{\alpha +z},aq^{\alpha -z};q)_{\infty }}{%
(aq^{z},aq^{-z};q)_{\infty }}\sum_{k=0}^{n}\frac{(q^{\alpha },abq^{\alpha
-1},acq^{\alpha -1},adq^{\alpha -1};q)_{k}}{(q,\gamma q^{2\alpha
-2},aq^{\alpha +z},aq^{\alpha -z};q)_{k}}q^{k},
\end{equation*}%
and $\gamma =abcd.$ (This is an analogue of the power series expansion; see
\cite{At:Su:DHF}, \cite[Exercises~2.9--2.11]{Su4}, and \cite{Su2} for
properties of the generalized powers.)

Apply the Askey--Wilson operator to $u_{n}^{\alpha }(x,y)$ to obtain
\begin{equation*}
(L_{0}(x)+\lambda _{\nu })u_{n}^{\alpha }(x,y)=\lambda _{\nu
}\sum_{m=0}^{n}c_{m}v_{m}[x(s)-x(\xi )]^{(\alpha +m)}
+\sum_{m=0}^{n}c_{m}v_{m}\ L_{0}(x)[x(s)-x(\xi )]^{(\alpha +m)},
\end{equation*}%
since $v_{m}$ is independent of $x.$ By \cite{At:Su:DHF}, we have
\begin{align}
L_{0}(x)[x(s)-x(\xi )]^{(\alpha +m)}
& =\gamma (\alpha +m)\gamma (\alpha +m-1)\sigma (\xi -\alpha
-m+1)[x(s)-x(\xi -1)]^{(\alpha +m-2)}  \notag \\
& \qquad +\gamma (\alpha +m)\tau _{\alpha +m-1}(\xi -\alpha -m+1)[x(s)-x(\xi
-1)]^{(\alpha +m-1)}  \notag \\
& \qquad -\lambda _{\alpha +m}[x(s)-x(\xi )]^{(\alpha +m)}.
\notag
\end{align}%
We use the same notations as in \cite{At:Su:DHF},
\cite[Exercise~2.25]{Su4},
or \cite{Su2}. Choose $a_{0}:=\xi -\alpha -m+1$ to be a
root of the equation $\sigma (a_{0})=0.$ Then $\xi =a_{0}+\alpha +m-1$, and
one obtains
\begin{align}
& (L_{0}(x)+\lambda _{\nu })u_{n}^{\alpha
}(x,y)=\sum_{m=0}^{n}c_{m}v_{m}\gamma (\alpha +m)\tau _{\alpha
+m-1}(a_{0})[x(s)-x(a_{0}+\alpha +m-2)]^{(\alpha +m-1)} \notag \\
& \qquad +\sum_{m=0}^{n}c_{m}v_{m}(\lambda _{\nu }-\lambda _{\alpha
+m})[x(s)-x(a_{0}+\alpha +m-1)]^{(\alpha +m)}  \notag \\
& =c_{0}v_{0}\gamma (\alpha )\tau _{\alpha -1}(a_{0})[x(s)-x(a_{0}+\alpha
-2)]^{(\alpha -1)}  \notag \\
& \qquad +\sum_{m=1}^{n}c_{m}v_{m}\gamma (\alpha +m)\tau _{\alpha
+m-1}(a_{0})[x(s)-x(a_{0}+\alpha +m-2)]^{(\alpha +m-1)}  \notag \\
& \qquad \qquad +\sum_{m=0}^{n}c_{m}v_{m}(\lambda _{\nu }-\lambda _{\alpha
+m})[x(s)-x(a_{0}+\alpha +m-1)]^{(\alpha +m)}.  \notag
\end{align}%
Letting $m=k+1,$ we get
\begin{multline}
(L_{0}(x)+\lambda _{\nu })u_{n}^{\alpha }(x,y)=c_{0}v_{0}\gamma (\alpha
)\tau _{\alpha -1}(a_{0})[x(s)-x(a_{0}+\alpha -2)]^{(\alpha -1)}
\label{L_0-eqn} \\
\kern4cm
+\sum_{k=0}^{n-1}c_{k+1}v_{k+1}\gamma (\alpha +k+1)\tau _{\alpha
+k}(a_{0})[x(s)-x(a_{0}+\alpha +k-1)]^{(\alpha +k)}  \\
+\sum_{k=0}^{n}c_{k}v_{k}(\lambda _{\nu }-\lambda _{\alpha
+k})[x(s)-x(a_{0}+\alpha +k-1)]^{(\alpha +k)}.
\end{multline}%
Note that for%
\begin{equation*}
v_{k}=\sum_{l=0}^{k}e_{l},\qquad e_{l}:=\frac{(aq^{\alpha +z},aq^{\alpha
-z};q)_{\infty }}{(aq^{z},aq^{-z};q)_{\infty }}\frac{(q^{\alpha
},abq^{\alpha -1},acq^{\alpha -1},adq^{\alpha -1};q)_{l}}{(q,\gamma
q^{2\alpha -2},aq^{\alpha +z},aq^{\alpha -z};q)_{l}}q^{l}
\end{equation*}%
one has
\begin{equation*}
v_{k+1}=v_{k}+e_{k+1}\text{ \ \ \ \ and \ \ \ \ }v_{0}=e_{0.}
\end{equation*}%
After choosing $\lambda _{\nu }=\lambda _{\alpha +n},$
equation~(\ref{L_0-eqn}) becomes
\begin{align*}
(L_{0}(x)+\lambda _{\alpha +n})u_{n}^{\alpha
}(x,y)&=\sum_{k=-1}^{n-1}c_{k+1}e_{k+1}\gamma (\alpha +k+1)\tau _{\alpha
+k}(a_{0})[x(s)-x(a_{0}+\alpha +k-1)]^{(\alpha +k)} \\
& \kern1cm
+\sum_{k=0}^{n-1}c_{k+1}v_{k}\gamma (\alpha +k+1)\tau _{\alpha
+k}(a_{0})[x(s)-x(a_{0}+\alpha +k-1)]^{(\alpha +k)} \\
& \kern1cm +\sum_{k=0}^{n-1}c_{k}v_{k}(\lambda _{\alpha +n}-\lambda
_{\alpha +k})[x(s)-x(a_{0}+\alpha +k-1)]^{(\alpha +k)}.
\end{align*}%
The latter two sums vanish if
\begin{equation*}
c_{k+1}\gamma (\alpha +k+1)\tau _{\alpha +k}(a_{0})=c_{k}(\lambda _{\alpha
+n}-\lambda _{\alpha +k}).
\end{equation*}%
Therefore,
\begin{multline*}
(L_{0}(x)+\lambda _{\alpha +n})u_{n}^{\alpha
}(x,y)=\sum_{k=-1}^{n-1}c_{k+1}e_{k+1}\gamma (\alpha +k+1)\tau _{\alpha
+k}(a_{0})[x(s)-x(a_{0}+\alpha +k-1)]^{(\alpha +k)} \\
=\sum_{m=0}^{n}c_{m}e_{m}\gamma (\alpha +m)\tau _{\alpha
+m-1}(a_{0})[x(s)-x(a_{0}+\alpha +m-2)]^{(\alpha +m-1)}=:f_{n}^{\alpha
}(x,y).
\end{multline*}%
Finally, we show that the function $f_{n}^{\alpha }(x,y)$ is, up to a factor,
the $n$-th ordinary Askey--Wilson polynomial. The generalized powers have
the property (see \cite{Su4})
\begin{equation*}
\lbrack x(s)-x(z)]^{(n+1)}=[x(s)-x(z)][x(s)-x(z-1)]^{(n)},
\end{equation*}%
which leads to
\begin{equation*}
f_{n}^{\alpha }(x,y)=\sum_{m=0}^{n}c_{m}e_{m}\gamma (\alpha +m)\tau _{\alpha
+m-1}(a_{0})\frac{[x(s)-x(a_{0}+\alpha +m-1)]^{(\alpha +m)}}{%
[x(s)-x(a_{0}+\alpha +m-1)]}.
\end{equation*}%
Moreover,
\begin{align*}
c_{m}[x(s)-x(a_{0}& +\alpha +m-1)]^{(\alpha +m)} \\
& =c_{0}\frac{(q^{-n},\gamma q^{2\alpha +n-1};q)_{m}}{(q^{\alpha
+1},abq^{\alpha },acq^{\alpha },adq^{\alpha };q)_{m}}q^{m}\
[x(s)-x(a_{0}+\alpha +m-1)]^{(\alpha +m)} \\
& =c_{0}\,\varphi _{m}(x)\,[x(s)-x(a_{0}+\alpha -1)]^{(\alpha )},
\end{align*}%
where, by definition,
\begin{equation*}
\varphi _{m}(x):=\frac{(aq^{s},aq^{-s};q)_{\infty }}{(aq^{\alpha
+s},aq^{\alpha -s};q)_{\infty }}\frac{(q^{-n},\gamma q^{2\alpha
+n-1},aq^{\alpha +s},aq^{\alpha -s};q)_{m}}{(q^{\alpha +1},abq^{\alpha
},acq^{\alpha },adq^{\alpha };q)_{m}}q^{m}.
\end{equation*}%
Therefore,
\begin{align}
f_{n}^{\alpha }(x,y)& =\frac{(aq^{s},aq^{-s};q)_{\infty }}{(aq^{\alpha
+s},aq^{\alpha -s};q)_{\infty }}\frac{(q^{-n},\gamma q^{2\alpha
+n-1},aq^{\alpha +s},aq^{\alpha -s};q)_{m}}{(q^{\alpha +1},abq^{\alpha
},acq^{\alpha },adq^{\alpha };q)_{m}}q^{m}  \notag \\
& \qquad \times \frac{(aq^{\alpha +z},aq^{\alpha -z};q)_{\infty }}{%
(aq^{z},aq^{-z};q)_{\infty }}\frac{(q^{\alpha },abq^{\alpha -1},acq^{\alpha
-1},adq^{\alpha -1};q)_{m}}{(q,\gamma q^{2\alpha -2},aq^{\alpha
+z},aq^{\alpha -z};q)_{m}}q^{m}  \notag \\
& \qquad \qquad \times \frac{\gamma (\alpha +m)\tau _{\alpha +m-1}(a_{0})}{%
\,[x(s)-x(a_{0}+\alpha +m-1)]}.
\label{FLemma1}
\end{align}%
Recall that \ $a=q^{a_{0}}$ and
\begin{align*}
\gamma (\alpha +m)& =q^{-\frac{\alpha +m-1}{2}}\,\frac{1-q^{\alpha +m}}{1-q},
\\
x(s)-x(a_{0}+\alpha +m-1)& =-\frac{1}{2a}q^{-\alpha -m+1}(1-aq^{\alpha
-s+m-1})(1-aq^{\alpha +s+m-1}),\text{ \ \ \ \ } \\
\tau _{\alpha +m-1}(a_{0})& =\frac{2}{a(1-q)}q^{-2(\alpha +m-1)+\frac{\alpha
+m}{2}}(1-abq^{\alpha +m-1})(1-acq^{\alpha +m-1})(1-adq^{\alpha +m-1}),
\end{align*}%
which allows us to simplify the last term of (\ref{FLemma1}) to
\begin{equation*}
q^{m}\frac{\gamma (\alpha +m)\tau _{\alpha +m-1}(a_{0})}{\,[x(s)-x(a_{0}+%
\alpha +m-1)]}=-4q^{\frac{3}{2}-\alpha }\frac{1-q^{\alpha +m}}{1-q}\frac{%
(1-abq^{\alpha +m-1})(1-acq^{\alpha +m-1})(1-adq^{\alpha +m-1})}{%
(1-aq^{\alpha -s+m-1})(1-aq^{\alpha +s+m-1})}.
\end{equation*}%
Thus $f_{n}^{\alpha }(x,y)$ becomes
\begin{align}
f_{n}^{\alpha }(x,y)& =\frac{-4q^{\frac{3}{2}-\alpha }}{(1-q)^{2}}\frac{%
(aq^{s},aq^{-s},aq^{\alpha +z},aq^{\alpha -z};q)_{\infty }}{(aq^{\alpha
+s-1},aq^{\alpha -s-1},aq^{z},aq^{-z};q)_{\infty }}\ (q^{\alpha
},abq^{\alpha -1},acq^{\alpha -1},adq^{\alpha -1};q)_{1} \\
& \qquad \times \sum_{m=0}^{n}\frac{(q^{-n},\gamma q^{2\alpha
+n-1},aq^{\alpha +s-1},aq^{\alpha -s-1};q)_{m}}{(q,\gamma q^{2\alpha
-2},aq^{\alpha +z},aq^{\alpha -z};q)_{m}}q^{m}  \notag \\
& =\frac{-4q^{\frac{3}{2}-\alpha }}{(1-q)^{2}}\frac{(aq^{s},aq^{-s},aq^{%
\alpha +z},aq^{\alpha -z};q)_{\infty }}{(aq^{\alpha +s-1},aq^{\alpha
-s-1},aq^{z},aq^{-z};q)_{\infty }}\medskip \ (q^{\alpha },abq^{\alpha
-1},acq^{\alpha -1},adq^{\alpha -1};q)_{1}\,  \notag \\
& \qquad \times p_{n}(x;aq^{\alpha -1},bcdq^{\alpha -1},q^{1+z},q^{1-z}),
\notag
\end{align}%
which completes the proof of the lemma.

\subsection{Proof of Lemma \protect\ref{lem:2}}

Consider the equation
\begin{equation*}
(L_{1}(y)+\lambda _{\nu })u_{n}^{\alpha }(x,y)=g_{n}^{\alpha }(x,y),
\end{equation*}%
and rewrite $u_{n}^{\alpha }(x,y)$ in the form
\begin{equation*}
u_{n}^{\alpha }(x,y)=\sum_{m=0}^{n}c_{m}^{\alpha }(aq^{\alpha +s},aq^{\alpha
-s};q)_{m}\frac{(aq^{s},aq^{-s};q)_{\infty }}{(aq^{\alpha +s},aq^{\alpha
-s};q)_{\infty }}v_{m}^{\alpha }(y), 
\end{equation*}%
where
\begin{equation*}
c_{m}^{\alpha }=\frac{(q^{-n},\gamma q^{2\alpha +n-1};q)_{m}}{(q^{\alpha
+1},abq^{\alpha },acq^{\alpha },adq^{\alpha };q)_{m}}q^{m},\text{ \ \ \ \ }%
\gamma =abcd,
\end{equation*}%
and
\begin{equation*}
v_{m}^{\alpha }(y)=\frac{(aq^{\alpha +z},aq^{\alpha -z};q)_{\infty }}{%
(aq^{z},aq^{-z};q)_{\infty }}\sum_{k=0}^{m}\frac{(q^{\alpha },abq^{\alpha
-1},acq^{\alpha -1},adq^{\alpha -1};q)_{k}}{(q,\gamma q^{2\alpha
-2},aq^{\alpha +z},aq^{\alpha -z};q)_{k}}q^{k}.
\end{equation*}%
Apply the Askey--Wilson operator $L_{1}(y):=L\left( y;q/a,q/b,q/c,q/d\right)
$ to $u_{n}^{\alpha }(x,y)$ to obtain
\begin{equation*}
(L_{1}(y)+\lambda _{\nu })u_{n}^{\alpha }(x,y)=\sum_{m=0}^{n}c_{m}^{\alpha
}(aq^{\alpha +s},aq^{\alpha -s};q)_{m}\frac{(aq^{s},aq^{-s};q)_{\infty }}{%
(aq^{\alpha +s},aq^{\alpha -s};q)_{\infty }}\ \left( L_{1}(y)+\lambda _{\nu
}\right) v_{m}^{\alpha }(y).
\end{equation*}

Let
\begin{equation*}
v_{m}^{\alpha }(y):=\sum_{k=0}^{m}\frac{c_{k}}{[x(s)-x(\xi )]^{(\alpha +k)}}
\end{equation*}%
in analogy with \cite{At:Su:DHF}. Then
\begin{equation*}
(L_{1}(y)+\lambda _{\nu })v_{m}^{\alpha }(y)=\lambda _{\nu }\sum_{k=0}^{m}%
\frac{c_{k}}{[x(s)-x(\xi )]^{(\alpha +k)}}+\sum_{k=0}^{m}c_{k}\
L_{1}(y)\left( \frac{1}{[x(s)-x(\xi )]^{(\alpha +k)}}\right) .
\end{equation*}%
By \cite{At:Su:DHF}, we have
\begin{multline}
L_{1}(y)\left( \frac{1}{[x(s)-x(\xi )]^{(\alpha +k)}}\right) =\frac{\gamma
(\alpha +k)\gamma (\alpha +k+1)\sigma (\xi +1)}{[x(z)-x(\xi +1)]^{(\alpha
+k+2)}}  \notag \\
-\frac{\gamma (\alpha +k)\tau _{-\alpha -k-1}(\xi +1)}{[x(z)-x(\xi
)]^{(\alpha +k+1)}}-\frac{\lambda _{-\alpha -k}}{[x(z)-x(\xi )]^{(\alpha +k)}%
}  \notag
\end{multline}%
(see also \cite[Exercise~2.25]{Su4}). Upon choosing $a_{0}:=\xi +1$ to
be a root of the equation $\sigma (a_{0})=0,$ we obtain
\begin{align}
(L_{1}(y)+\lambda _{\nu })v_{m}^{\alpha }(y)& =\lambda _{\nu }\sum_{k=0}^{m}%
\frac{c_{k}}{[x(s)-x(a_{0})]^{(\alpha +k)}} \notag \\
& \qquad -\sum_{k=0}^{m}c_{k}\left( \frac{\gamma (\alpha +k)\tau _{-\alpha
-k-1}(a_{0})}{[x(z)-x(a_{0}-1)]^{(\alpha +k+1)}}+\frac{\lambda _{-\alpha -k}%
}{[x(z)-x(a_{0}-1)]^{(\alpha +k)}}\right)   \notag \\
& =\sum_{k=0}^{m}\frac{c_{k}\left( \lambda _{\nu }-\lambda _{-\alpha
-k}\right) }{[x(z)-x(a_{0}-1)]^{(\alpha +k)}}-\sum_{k=0}^{m}\frac{%
c_{k}\,\gamma (\alpha +k)\tau _{-\alpha -k-1}(a_{0})}{[x(z)-x(a_{0}-1)]^{(%
\alpha +k+1)}}  \notag \\
& =\frac{c_{0}\left( \lambda _{\nu }-\lambda _{-\alpha }\right) }{%
[x(z)-x(a_{0}-1)]^{(\alpha )}}+\sum_{k=1}^{m}\frac{c_{k}\left( \lambda _{\nu
}-\lambda _{-\alpha -k}\right) }{[x(z)-x(a_{0}-1)]^{(\alpha +k)}}  \notag \\
& \qquad -\frac{c_{m}\,\gamma (\alpha +m)\tau _{-\alpha -m-1}(a_{0})}{%
[x(z)-x(a_{0}-1)]^{(\alpha +m+1)}}-\sum_{k=0}^{m-1}\frac{c_{k}\,\gamma
(\alpha +k)\tau _{-\alpha -k-1}(a_{0})}{[x(z)-x(a_{0}-1)]^{(\alpha +k+1)}}.
\notag
\end{align}%
Now choose $\lambda _{\nu }=\lambda _{-\alpha }$ and let $k=l+1.$
Then we obtain%
\begin{multline*}
(L_{1}(y)+\lambda _{\nu })v_{m}^{\alpha }(y)=-\frac{c_{m}\,\gamma
(\alpha +m)\tau _{-\alpha -m-1}(a_{0})}{[x(z)-x(a_{0}-1)]^{(\alpha +m+1)}} \\
+\sum_{l=0}^{m-1}\frac{c_{l+1}\left( \lambda _{-\alpha }-\lambda
_{-\alpha -l-1}\right) }{[x(z)-x(a_{0}-1)]^{(\alpha +l+1)}}-\sum_{l=0}^{m-1}%
\frac{c_{l}\,\gamma (\alpha +l)\tau _{-\alpha -l-1}(a_{0})}{%
[x(z)-x(a_{0}-1)]^{(\alpha +l+1)}}.
\end{multline*}%
The latter two sums vanish if
\begin{equation*}
c_{l+1}\left( \lambda _{-\alpha }-\lambda _{-\alpha -l-1}\right)
=c_{l}\,\gamma (\alpha +l)\tau _{-\alpha -l-1}(a_{0}).
\end{equation*}%
In that case, we have
\begin{equation*}
(L_{1}(y)+\lambda _{\nu })v_{m}^{\alpha }(y)=-\frac{c_{m}\,\gamma (\alpha
+m)\tau _{-\alpha -m-1}(a_{0})}{[x(z)-x(a_{0}-1)]^{(\alpha +m+1)}}=-\frac{%
c_{m+1}\left( \lambda _{-\alpha }-\lambda _{-\alpha -m-1}\right) }{%
[x(z)-x(a_{0}-1)]^{(\alpha +m+1)}}=:h_{m}^{\alpha }(y).
\end{equation*}%
Here,
\begin{align*}
\frac{c_{m+1}}{[x(z)-x(a_{0}-1)]^{(\alpha +m+1)}}& =\frac{c_{0}}{%
[x(z)-x(a_{0}-1)]^{(\alpha )}}\varphi _{m+1}(z), \\
\varphi _{m+1}(z)& =\frac{(q^{\alpha },abq^{\alpha -1},acq^{\alpha
-1},adq^{\alpha -1};q)_{m+1}}{(q,\gamma q^{2\alpha -2},aq^{\alpha
+z},aq^{\alpha -z};q)_{m+1}}q^{m+1}, \\
\frac{c_{0}}{[x(z)-x(a_{0}-1)]^{(\alpha )}}& =\frac{(aq^{\alpha
+z},aq^{\alpha -z};q)_{\infty }}{(aq^{z},aq^{-z};q)_{\infty }}, \\
\lambda _{-\alpha }-\lambda _{-\alpha -m-1}& =\frac{4}{(1-q)^{2}\gamma }q^{%
\frac{7}{2}-\alpha -m}(1-q^{m+1})(1-\gamma q^{2\alpha +m-2})
\end{align*}%
and
\begin{equation*}
h_{m}^{\alpha }(y)=-\frac{4q^{\frac{9}{2}-\alpha }}{(1-q)^{2}\gamma }\frac{%
(aq^{\alpha +z},aq^{\alpha -z};q)_{\infty }}{(aq^{z},aq^{-z};q)_{\infty }}%
\frac{(q^{\alpha },abq^{\alpha -1},acq^{\alpha -1},adq^{\alpha -1};q)_{m+1}}{%
(q,\gamma q^{2\alpha -2};q)_{m}(aq^{\alpha +z},aq^{\alpha -z};q)_{m+1}}.
\end{equation*}%
Therefore,
\begin{align}
(L_{1}(y)& +\lambda _{\nu })u_{n}^{\alpha }(x,y) \notag \\
& =\sum_{m=0}^{n}c_{m}^{\alpha }(aq^{\alpha +s},aq^{\alpha -s};q)_{m}\frac{%
(aq^{s},aq^{-s};q)_{\infty }}{(aq^{\alpha +s},aq^{\alpha -s};q)_{\infty }}%
L_{1}(y)v_{m}^{\alpha }(y)  \notag \\
& =-\sum_{m=0}^{n}c_{m}^{\alpha }(aq^{\alpha +s},aq^{\alpha -s};q)_{m}\frac{%
(aq^{s},aq^{-s};q)_{\infty }}{(aq^{\alpha +s},aq^{\alpha -s};q)_{\infty }}
\notag \\
& \qquad \times \frac{4q^{\frac{9}{2}-\alpha }}{(1-q)^{2}\gamma }
\frac{%
(aq^{\alpha +z},aq^{\alpha -z};q)_{\infty }}{(aq^{z},aq^{-z};q)_{\infty }}%
\frac{(q^{\alpha },abq^{\alpha -1},acq^{\alpha -1},adq^{\alpha -1};q)_{m+1}}{%
(q,\gamma q^{2\alpha -2};q)_{m}(aq^{\alpha +z},aq^{\alpha -z};q)_{m+1}}
\notag \\
& =-\frac{4q^{\frac{9}{2}-\alpha }}{(1-q)^{2}\gamma }\frac{%
(aq^{s},aq^{-s},aq^{\alpha +z+1},aq^{\alpha -z+1};q)_{\infty }}{(aq^{\alpha
+s},aq^{\alpha -s},aq^{z},aq^{-z};q)_{\infty }}  \notag \\
& \qquad \times \sum_{m=0}^{n}\frac{(q^{-n},\gamma q^{2\alpha
+n-1},aq^{s},aq^{-s};q)_{m}}{(q^{\alpha +1},abq^{\alpha },acq^{\alpha
},adq^{\alpha };q)_{m}}q^{m}  \notag \\
& \qquad \times \frac{(1-q^{\alpha })(1-abq^{\alpha -1})(1-acq^{\alpha
-1})(1-adq^{\alpha -1})}{(q,\gamma q^{2\alpha -2};q)_{m}(aq^{\alpha
+z+1},aq^{\alpha -z+1};q)_{m}}  \notag \\
& =-\frac{4q^{\frac{9}{2}-\alpha }}{(1-q)^{2}\gamma }\frac{%
(aq^{s},aq^{-s},aq^{\alpha +z+1},aq^{\alpha -z+1};q)_{\infty }}{(aq^{\alpha
+s},aq^{\alpha -s},aq^{z},aq^{-z};q)_{\infty }}  \notag \\
& \qquad \times \left( q^{\alpha },abq^{\alpha -1},acq^{\alpha
-1},adq^{\alpha -1};q\right) _{1}  \notag \\
& \qquad \times \text{ }_{4}\varphi _{3}\!\left( \!\!%
\begin{array}{c}
q^{-n},\gamma q^{2\alpha +n-1},aq^{\alpha +s},aq^{\alpha -s} \\[0.1cm]
\gamma q^{2\alpha -2},aq^{\alpha +z+1},aq^{\alpha -z+1}%
\end{array}%
\!\!;q,\,q\!\right)   \notag \\
& =-\frac{4q^{\frac{9}{2}-\alpha }}{(1-q)^{2}\gamma }\frac{%
(aq^{s},aq^{-s},aq^{\alpha +z+1},aq^{\alpha -z+1};q)_{\infty }}{(aq^{\alpha
+s},aq^{\alpha -s},aq^{z},aq^{-z};q)_{\infty }}  \notag \\
& \qquad \times \left( q^{\alpha },abq^{\alpha -1},acq^{\alpha
-1},adq^{\alpha -1};q\right) _{1}\times \text{ }p_{n}(x;aq^{\alpha
},bcdq^{\alpha -2},q^{1+z},q^{1-z}).  \notag
\end{align}%
This completes the proof of the lemma.

\subsection{Proof of Lemma \protect\ref{lem:3}}

The structure of the Askey--Wilson operator in (\ref{in9}) and the basic
hypergeometric series representation (\ref{AskeyWilsonPlinomials}) suggest
to look for a $4$-term relation of the form%
\begin{multline}
K_{1}\ {}_{4}\varphi _{3}\!\left( \!\!%
\begin{array}{c}
A,B,C,D \\[0.1cm]
F,G,H%
\end{array}%
\!\!;q,\,q\!\right) +K_{2}\ {}_{4}\varphi _{3}\!\left( \!\!%
\begin{array}{c}
A,B,Cq,D/q \\[0.1cm]
F,G,H%
\end{array}%
\!\!;q,\,q\!\right)    \\
+K_{3}\ {}_{4}\varphi _{3}\!\left( \!\!%
\begin{array}{c}
A,B,C/q,Dq \\[0.1cm]
F,G,H%
\end{array}%
\!\!;q,\,q\!\right) +K_{4}\ {}_{4}\varphi _{3}\!\left( \!\!%
\begin{array}{c}
A,B,C/q,D/q \\[0.1cm]
F,G/q,H/q%
\end{array}%
\!\!;q,\,q\!\right) =0,
\label{4TermRecurrenceAWPols}
\end{multline}%
for some undetermined coefficients $K_{1},$ $K_{2},$ $K_{3}$ and $K_{4}$
(up to a common factor). Doing a term-wise comparison,
we may hope to find $K_{1},$ $K_{2},$ $K_{3},$ $K_{4}$
which satisfy%
\begin{align}
K_{1}(1-C)& (1-D)(1-Cq^{k-1})(1-Dq^{k-1})(1-G/q)(1-H/q)  \notag \\
& +K_{2}(1-Cq^{k})(1-Cq^{k-1})(1-D/q)(1-D)(1-G/q)(1-H/q)  \notag \\
& +K_{3}(1-Dq^{k})(1-Dq^{k-1})(1-C/q)(1-C)(1-G/q)(1-H/q)  \notag \\
& +K_{4}(1-C/q)(1-C)(1-D/q)(1-D)(1-Gq^{k-1})(1-Hq^{k-1})=0.  \notag
\end{align}%
If we are successful, then the above equation does indeed imply
the contiguous relation \eqref{4TermRecurrenceAWPols}.
In the equation, we compare coefficients of powers of $q^{k}$.
This yields a system of 3 linear equations in the 4 unknowns
$K_{1},$ $K_{2},$ $K_{3},$ $K_{4}$.
With the help of \textsl{Mathematica}, we obtain the solution
\begin{align*}
K_{1}& =\frac{(C-q)(D-q)(-GH-CDq+CGq+DGq+CHq+DHq-GHq-CDq^{2})}{%
(G-q)(H-q)(Cq-D)(Dq-C)}, \\
K_{2}& =\frac{(C-1)(D-G)(D-H)(C-q)q}{(D-C)(G-q)(H-q)(Cq-D)}, \\
K_{3}& =\frac{(D-1)(C-G)(C-H)(D-q)q}{(C-D)(G-q)(H-q)(Dq-C)},
\end{align*}%
where the free parameter $K_{4}$ was chosen to be $1$ (see Appendix~A for
the \textsl{Mathematica} code). The required $4$-term contiguous relation is
then given by
\begin{align}
& \frac{(C-q)(D-q)(-GH-CDq+CGq+DGq+CHq+DHq-GHq-CDq^{2})}{%
(G-q)(H-q)(Cq-D)(Dq-C)}\text{ }  \notag \\
& \qquad \quad \quad \quad  \times \text{ }_{4}\varphi _{3}\!\left( \!\!%
\begin{array}{c}
A,B,C,D \\[0.1cm]
F,G,H%
\end{array}%
\!\!;q,\,q\!\right) +{}_{4}\varphi _{3}\!\left( \!\!%
\begin{array}{c}
A,B,C/q,D/q \\[0.1cm]
F,G/q,H/q%
\end{array}%
\!\!;q,\,q\!\right)   \notag \\
& \qquad \quad +\frac{(C-1)(D-G)(D-H)(C-q)q}{(D-C)(G-q)(H-q)(Cq-D)}\text{ }%
_{4}\varphi _{3}\!\left( \!\!%
\begin{array}{c}
A,B,Cq,D/q \\[0.1cm]
F,G,H%
\end{array}%
\!\!;q,\,q\!\right)   \notag \\
& \qquad \quad +\frac{(D-1)(C-G)(C-H)(D-q)q}{(C-D)(G-q)(H-q)(Dq-C)}%
\text{ }_{4}\varphi _{3}\!\left( \!\!%
\begin{array}{c}
A,B,C/q,Dq \\[0.1cm]
F,G,H%
\end{array}%
\!\!;q,\,q\!\right) =0.
\label{4TermsRecurrenceAWSolved}
\end{align}%
(This $4$-term contiguous relation for the $_{4}\varphi _{3}$-functions can
be extended to an arbitrary $_{r}\psi _{s}$-function, see Appendix~A for
more details.)

When $qABCD=FGH$, in view of the structure of the Askey--Wilson operator in (%
\ref{in9}), equation~(\ref{4TermsRecurrenceAWSolved}) should become%
\begin{align*}
& (L(x)+\lambda )\ _{4}\varphi _{3}\!\left( \!\!%
\begin{array}{c}
A,B,Cq,D/q \\[0.1cm]
F,G,H%
\end{array}%
\!\!;q,\,q\!\right)  \\
& =\frac{\sigma (-s)}{\Delta x(s)\nabla x_{1}(s)}\ {}_{4}\varphi
_{3}\!\left( \!\!%
\begin{array}{c}
A,B,Cq,D/q \\[0.1cm]
F,G,H%
\end{array}%
\!\!;q,\,q\!\right) +\frac{\sigma (s)}{\nabla x(s)\nabla x_{1}(s)}{}\
_{4}\varphi _{3}\!\left( \!\!%
\begin{array}{c}
A,B,C/q,Dq \\[0.1cm]
F,G,H%
\end{array}%
\!\!;q,\,q\!\right)  \\
& \qquad +\frac{\lambda \Delta x(s)\nabla x(s)\nabla x_{1}(s)-\sigma
(s)\Delta x(s)-\sigma (-s)\nabla x(s)}{\Delta x(s)\nabla x(s)\nabla x_{1}(s)%
}{}\ _{4}\varphi _{3}\!\left( \!\!%
\begin{array}{c}
A,B,C,D \\[0.1cm]
F,G,H%
\end{array}%
\!\!;q,\,q\!\right) .
\end{align*}%
Equating coefficients, one obtains
\begin{multline*}
(1-C/q)(1-D/q)(D-C)(GH+CDq-CGq-DGq-CHq-DHq+GHq+CDq^{2}) \\
=\frac{2qa^{3}}{1-q}\left( \sigma (-s)\nabla x(s)+\sigma (s)\Delta
x(s)-\lambda \Delta x(s)\nabla x(s)\nabla x_{1}(s)\right)
\end{multline*}%
and%
\begin{align*}
(D-C)(D-C/q)(-C+D/q)&=\frac{2aq^{1/2}}{1-q}\nabla x_{1}(s)\frac{2a}{1-q}%
\Delta x(s)\frac{2a}{1-q}\nabla x(s), \\
(C-1)(G-D)(H-D)(q-C)&=-qa^{2}\sigma (-s), \\
(D-C)(-C+D/q)&=\frac{2aq^{1/2}}{1-q}\nabla x_{1}(s)\frac{2a}{1-q}\Delta
x(s), \\
(D-1)(G-C)(H-C)(q-D)&=-qa^{2}\sigma (s), \\
(D-C)(D-C/q)&=\frac{2aq^{1/2}}{1-q}\nabla x_{1}(s)\frac{2a}{1-q}\nabla x(s),
\\
(G-q)(H-q)&=q^{2}\frac{(1-q)^{2}}{4q^{3/2}}\lambda .
\end{align*}%
This gives the required formula (\ref{DiffDiffFormulaAW}) for the
Askey--Wilson operator with
\begin{equation*}
\sigma (s)=q^{-2s}\left( q^{s}-a\right) \left( q^{s}-a/q\right) \left(
q^{s}-c\right) \left( q^{s}-d\right) ,\quad \lambda =\frac{4q^{3/2}}{%
(1-q)^{2}}\left( 1-ac/q\right) \left( 1-ad/q\right) .
\end{equation*}%
The proof of the lemma is complete.

\appendix

\section{4-Term Contiguous Relations}

In order to derive the contiguous relation \eqref{4TermsRecurrenceAWSolved},
one can use the following \textsl{Mathematica} program:\footnote{A
corresponding \textsl{Mathematica} notebook is
available on the article's website\newline
\texttt{http://www.mat.univie.ac.at/\lower0.5ex\hbox{\~{}}kratt/artikel/AssAWPols.html}.}

\smallskip
\MATH
\goodbreakpoint%
In[1]:= X1 = K1*(1 - C) (1 - D) (1 - C*K/q) (1 - D*K/q) (1 - G/q) (1 - H/q) +
\ \ \ \ \ \ \ \ \     K2*(1 - C*K) (1 - C*K/q) (1 - D/q) (1 - D) (1 - G/q) (1 - H/q) +
\ \ \ \ \ \ \ \ \     K3*(1 - D*K) (1 - D*K/q) (1 - C/q) (1 - C) (1 - G/q) (1 - H/q) +
\ \ \ \ \ \ \ \ \     K4*(1 - C/q) (1 - C) (1 - D/q) (1 - D) (1 - G*K/q) (1 - H*K/q) ;
\ \ \ \ \    X1 = Table[Coefficient[X1, K, i] == 0, {i, 0, 2}];
\ \ \ \ \    X1 = Solve[X1, {K1, K2, K3, K4}];
\ \ \ \ \    X1 = \MATHlbrace %
K1 -> Factor[K1/.X1[[1]]], K2 -> Factor[K2/.X1[[1]]],
\              K3 -> Factor[K3/.X1[[1]]], K4 -> Factor[K4/.X1[[1]]]%
\MATHrbrace %
\medskip
Out[1]= %
\MATHlbrace %
K1 \MATHhStrich \MATHgroesser %
(K4 (C - q) (D - q)
\goodbreakpoint%
\MATHgroesser           %
(G H + C D q - C G q - D G q - C H q - D H q + G H q + C D q\MATHhoch 2)) /
\MATHgroesser           %
((G - q) (H - q) (-D + C q) (C - D q)),
\goodbreakpoint%
\MATHgroesser        %
K2 \MATHhStrich \MATHgroesser %
- ((-1 + C) (D - G) (D - H) K4 (C - q) q) /
\goodbreakpoint%
\MATHgroesser           %
((C - D) (G - q) (H - q) (-D + C q)),
\goodbreakpoint%
\MATHgroesser        %
K3 \MATHhStrich \MATHgroesser %
- ((-1 + D) (C - G) (C - H) K4 (D - q) q) /
\goodbreakpoint%
\MATHgroesser           %
((C - D) (G - q) (H - q) (C - D q)),
\goodbreakpoint%
\MATHgroesser        %
K4 \MATHhStrich \MATHgroesser %
K4\MATHrbrace %
\goodbreakpoint%
\endMATH

\medskip
It is evident from the proof of \eqref{4TermsRecurrenceAWSolved}
that, actually, an extension for bilateral series
(see \cite[equation~(5.1.1)]{Ga:Ra} for the definition) with an arbitrary
number of parameters holds, namely:
\begin{align} \notag
&\frac{\left( c-q\right) \left( d-q\right) \left(
-gh-cdq+cgq+dgq+chq+dhq-ghq-cdq^{2}\right) }{\left( g-q\right) \left(
h-q\right) \left( cq-d\right) \left( dq-c\right) } \\
&\quad \quad \quad \qquad \times \ _{r}\psi _{s}\!\left( \!\!%
\begin{array}{c}
a_{1},\dots,a_{i},\ c,d \\[0.1cm]
b_{0},\dots,b_{k},\ g,h%
\end{array}%
\!\!;q,\,t\!\right)   \notag \\
&\quad \quad
+\frac{\left( c-1\right) \left( d-g\right) \left( d-h\right) \left(
c-q\right) q}{\left( d-c\right) \left( g-q\right) \left( h-q\right) \left(
cq-d\right) }\ \ _{r}\psi _{s}\!\left( \!\!%
\begin{array}{c}
a_{1},\dots,a_{i},\ cq,d/q \\[0.1cm]
b_{0},\dots,b_{k},\ g,h%
\end{array}%
\!\!;q,\,t\!\right)   \notag \\
&\quad \quad
+\frac{\left( d-1\right) \left( c-g\right) \left( c-h\right) \left(
d-q\right) q}{\left( c-d\right) \left( g-q\right) \left( h-q\right) \left(
dq-c\right) }\ \ _{r}\psi _{s}\!\left( \!\!%
\begin{array}{c}
a_{1},\dots,a_{i},\ c/q,dq \\[0.1cm]
b_{0},\dots,b_{k},\ g,h%
\end{array}%
\!\!;q,\,t\!\right)   \notag \\
&\qquad +\ _{r}\psi _{s}\!\left( \!\!%
\begin{array}{c}
a_{1},\dots,a_{i},\ c/q,d/q \\[0.1cm]
b_{0},\dots,b_{k},\ g/q,h/q%
\end{array}%
\!\!;q,\,t\!\right) =0.
\label{4Polynomials}
\end{align}%
Furthermore, in the same way, the following variation can be 
obtained:\footnote{Again, a
corresponding {\sl Mathematica} notebook is
available on the article's website\newline
{\tt http://www.mat.univie.ac.at/\lower0.5ex\hbox{\~{}}kratt/artikel/AssAWPols.html}.}
\begin{align} \notag
&\frac{\left( g-1\right) \left( h-1\right) \left(
-gh-cdq+cgq+dgq+chq+dhq-ghq-cdq^{2}\right) }{\left( c-1\right) \left(
d-1\right) \left( gq-h\right) \left( hq-g\right) } \\
&\quad \quad \quad \qquad \times \ _{r}\psi _{s}\!\left( \!\!%
\begin{array}{c}
a_{1},\dots,a_{i},\ c,d \\[0.1cm]
b_{0},\dots,b_{k},\ g,h%
\end{array}%
\!\!;q,\,t\!\right)   \notag \\
&\quad \quad
+\frac{\left( c-g\right) \left( d-g\right) \left( h-1\right) \left(
h-q\right) }{\left( c-1\right) \left( d-1\right) \left( h-g\right) \left(
gq-h\right) }\ \ _{r}\psi _{s}\!\left( \!\!%
\begin{array}{c}
a_{1},\dots,a_{i},\ c,d \\[0.1cm]
b_{0},\dots,b_{k},\ gq,h/q%
\end{array}%
\!\!;q,\,t\!\right)   \notag \\
&\quad \quad
+\frac{\left( c-h\right) \left( d-h\right) \left( g-1\right) \left(
g-q\right) }{\left( c-1\right) \left( d-1\right) \left( g-h\right) \left(
hq-g\right) }\ \ _{r}\psi _{s}\!\left( \!\!%
\begin{array}{c}
a_{1},\dots,a_{i},\ c,d \\[0.1cm]
b_{0},\dots,b_{k},\ g/q,hq%
\end{array}%
\!\!;q,\,t\!\right)   \notag \\
&\qquad +\ _{r}\psi _{s}\!\left( \!\!%
\begin{array}{c}
a_{1},\dots,a_{i},\ cq,dq \\[0.1cm]
b_{0},\dots,b_{k},\ gq,hq%
\end{array}%
\!\!;q,\,t\!\right) =0.
\label{4RationalF}
\end{align}

\section{An Inverse of the Askey--Wilson Operator}

The Askey--Wilson divided difference operator on the left-hand side of
equation~(\ref{DiffDiffFormulaAW}) can be inverted by the method of Ref.~%
\cite{As:Rah:Sus}. The end result is%
\begin{equation}
\frac{\left( q,q^{2};q\right) _{\infty }}{2\pi }\int_{-1}^{1}L\left(
x,y\right) \ p_{n}\left( x;a,b,c,d\right) \rho \left( x;a,b,c,d\right) \
dx=p_{n}\left( x;aq,b/q,c,d\right) ,
\end{equation}%
where $\rho \left( x;a,b,c,d\right) $ is the weight function of the
Askey--Wilson polynomials (\ref{AskeyWilsonPlinomials}) and the kernel is
given by%
\begin{multline}
L\left( x,y\right)  =\left( ac,ad,qce^{i\varphi },qde^{-i\varphi
};q\right) _{1}\ \frac{\left( be^{i\theta },be^{-i\theta },qde^{i\theta
},qde^{-i\theta },qae^{i\varphi },qae^{-i\varphi },qce^{i\varphi
},qce^{-i\varphi };q\right) _{\infty }}{\left( qe^{i\theta +i\varphi
},qe^{i\theta -i\varphi },qe^{i\varphi -i\theta },qe^{-i\theta -i\varphi
};q\right) _{\infty }}  \notag \\
\times \ _{8}\varphi _{7}\!\left( \!\!%
\begin{array}{c}
qde^{-i\varphi },q\sqrt{qde^{-i\varphi }},-q\sqrt{qde^{-i\varphi }}%
,qe^{i\theta -i\varphi },qe^{-i\theta -i\varphi },qd/c,q \\[0.1cm]
\sqrt{qde^{-i\varphi }},\sqrt{qde^{-i\varphi }},qde^{-i\theta },qde^{i\theta
},q^{2},qce^{-i\varphi },qde^{-i\varphi }%
\end{array}%
\!\!;q,\,ce^{i\varphi }\right) .
\end{multline}%
Here, $x=\cos \theta $ and $y=\cos \varphi .$ Computational details are left
to the reader.

\noindent \noindent \textbf{Acknowledgment.\/} We thank Mizan Rahman
for valuable discussions and encouragement.

\end{document}